\documentclass[11pt]{article}

\usepackage{amsmath, epsfig, cite,colordvi,xcolor}
\usepackage{amsthm}
\usepackage{amssymb}
\usepackage{amsfonts}
\usepackage{latexsym}
\usepackage{graphicx}
\usepackage{tikz}
\usepackage{float}
\usepackage{mathrsfs}[12pt]
\usepackage{bm, amscd,  amsfonts, amssymb, cite,texdraw}
\usepackage[colorlinks,urlcolor=purple,linkcolor=red,anchorcolor=green,citecolor=blue]{hyperref}

\newtheorem{thm}{Theorem}[section]

\newtheorem{lem}[thm]{Lemma}

\theoremstyle{definition}

\renewcommand{\qed}{{\hfill\rule{4pt}{7pt}}}
\def\pf{\noindent {\it Proof.} }
\def\gen{{\rm Gen}}

\def\qed{\hfill \rule{7pt}{7pt}}

\numberwithin{equation}{section}

\makeatletter \@addtoreset{equation}{section} \makeatother

\tikzstyle{every node}=[circle,inner sep=1pt,fill=white!60]
\tikzstyle{tn}=[shape=circle, draw, color=black!70]

\begin{document}
\begin{center}
{\bf \large A Context-free Grammar for

 Peaks and Double Descents of Permutations }

\vskip 9pt

Amy M. Fu

School of Mathematics\\
Shanghai University of Finance and Economics\\
 Shanghai 200433, P.R. China

 fu.mei@mail.shufe.edu.cn

\end{center}

\noindent{\bf Abstract.}
This paper is concerned with  the joint distribution of the number of exterior
peaks and the number of proper double descents over
permutations on $[n]=\{1,2,\ldots,n\}$.
The notion of exterior peaks of a permutation was introduced by
Aguiar, Bergeron and Nyman in their study of the peak algebra.
Gessel obtained the generating function of
the number of permutations on $[n]$ with a given number of exterior peaks.
On the other hand, by establishing   differential equations,
 Elizalde and Noy derived the generating function for the number of permutations on $[n]$ with a given number of
proper double descents.
Barry and Basset independently deduced the generating function of the number of permutations on $[n]$ with no proper double descents.
We find a context-free grammar which can be used to
compute the number of permutations on $[n]$ with a given number of exterior peaks and a given number of proper double descents. Based on the grammar,
the recurrence relation of the number of permutations on $[n]$ with a give number of exterior peaks can be easily obtained. Moreover, we use the grammatical calculus to derive
the generating function without solving  differential
equations. Our formula reduces
to the formulas of Gessel, Elizalde-Noy, Barry, and Basset.
Finally, from the grammar we establish a relationship between
our generating function and the generating function of the
joint distribution of the number of peaks and the number of
double descent derived by Carlitz and Scoville.



\section{Introduction}

The objective of this paper is to present a
grammatical approach to the joint distribution of  exterior peaks and proper double descents of permutations on $[n]$.
The notion of  exterior peaks was introduced by Aguiar, Bergeron and Nyman \cite{Aguiar-2004}.
Given a permutation $\pi=\pi_1\pi_2\cdots\pi_n$ on $[n]$,
an index $i$ is called an {\it exterior peak} if $\pi_1 > \pi_2$ for $i = 1$ or $\pi_{i-1}< \pi_i > \pi_{i+1}$
for $1 < i < n$. Let $T(n,k)$ be the number of  permutations on $[n]$ with $k$ exterior peaks and let
\begin{equation}\label{thenumberofpeaks}
  T_n(x)=\sum_{k\geq 0}T(n,k)x^k.
\end{equation}
 Gessel\cite{Gessel-OEIS} obtained the generating function
 of $T_n(x)$.
 \begin{thm}[Gessel\cite{Gessel-OEIS}]\label{lem1}
We have
\begin{align}\label{gf2}
\sum_{n=0}^\infty \frac{T_n(x)t^n}{n!}
=\frac{\sqrt{1-x}}{\sqrt{1-x}\cosh{(t\sqrt{1-x})}-\sinh{(t\sqrt{1-x}})}.
\end{align}
\end{thm}

The number of proper double descents  is a classical statistic on permutations, which has been extensively studied. An index $i$ of a permutation $\pi=\pi_1\pi_2\ldots\pi_n$ on $[n]$ is called  a {\it proper double descent} if $3\leq i\leq n$ and
$\pi_{i-2} > \pi_{i-1} > \pi_i$.
Denote $U(n,k)$
 the number of permutations on $[n]$ with $k$ proper double descents and let
 \begin{equation}\label{yny}
 U_n(y)=\sum_{k\geq 0}U(n,k)y^k.
 \end{equation}
 By establishing the following ordinary differential equations,
 \begin{align*}
  f''+(1-y)\left(
  f'+f\right)=0
\end{align*}
with $f(0)=1$ and $f'(0)=-1$, Elizalde and Noy\cite{Elizalde-2003} derived the generating function of $U_n(y)$.

 \begin{thm}[Elizalde and Noy \cite{Elizalde-2003}]\label{lem2}
We have
\begin{align}\label{gf3}
\sum_{n=0}^\infty \frac{U_n(y)t^n}{n!}
=\frac{2\sqrt{(y-1)(y+3)}e^{t/2\cdot(1-y+\sqrt{(y-1)(y+3)})}}{1+y+\sqrt{(y-1)(y+3)}-(1+y-\sqrt{(y-1)(y+3)})e^{t\sqrt{(y-1)(y+3)}}}.
\end{align}
\end{thm}
Barry\cite{Barry-2014} and Basset\cite{Basset-2014} also independently studied the generating function for the number of permutations on $[n]$ with no proper double descents.

\begin{thm}[Elizalde and Noy\cite{Elizalde-2003}, Barry \cite{Barry-2014}, Basset \cite{Basset-2014}]\label{lem3}
We have
  \begin{equation}\label{No-doubledescents-formula}
  \sum_{n=0}^{\infty}U(n,0)\frac{t^n}{n!}
    =\frac{\sqrt{3}}{2}\frac{e^{t/2}}{\cos{(\sqrt{3}t/2+\pi/6)}}.
  \end{equation}
\end{thm}
Note that Theorem \ref{lem3}  gives the explicit form of the following generating function
\begin{equation}\label{GFofnodoublefalls1}
  \sum_{n=0}^\infty U(n,0)\frac{t^n}{n!}=\left(\sum_{j=0}^\infty\cfrac{t^{3j}}{(3j)!}-\sum_{j=0}^\infty\cfrac{t^{3j+1}}{(3j+1)!}\right)^{-1},
\end{equation}
which appears in  \cite[pp.156--157]{David-1962}, \cite[Example 3, pp. 51]{Gessel-1977}, \cite[Exercise 5.2.17]{Goulden-2004} and \cite[pp. 126 and 260]{Stanley-1997}.
The right hand sides of \eqref{No-doubledescents-formula} and \eqref{GFofnodoublefalls1}  are equal, which can be seen from the following formula
\begin{equation*}
  \begin{split}
   \sum_{j=0}^\infty\cfrac{t^{3j}}{(3j)!}-\sum_{j=0}^\infty\cfrac{t^{3j+1}}{(3j+1)!}&=\frac{1}{3}\left(e^x+e^{\omega x}+e^{\omega^2x}\right)-\frac{1}{3}\left(e^x+\omega^2e^{\omega x}+\omega e^{\omega^2x}\right)\\[6pt]
   &=\frac{1}{3}(1-\omega^2)e^{\omega x}+\frac{1}{3}(1-\omega)e^{\omega^2 x},
  \end{split}
\end{equation*}
where $\omega=e^{2\pi\mathrm{i}/3}$.

To consider the joint distribution of the
number of exterior peaks and the number of
proper double descents over permutations, we
define $P_n(i, j)$ to be the number of   permutations on $[n]$ with
  $i$ exterior peaks and
$j$ proper double descents. Let
\begin{equation}\label{Pn}
P_n(x,y)=\sum_{i,j}
P_n(i, j)x^i y^j,
\end{equation}
where $0\leq j \leq n-1$ and  $2i+j\leq n$.
we find a context free grammar and a grammatical labeling of
permutations to generate the polynomials $P_n(x,y)$.
In fact, for the reason of the grammar, we need to define
 the polynomials {$P_n(x,y,z,w)$} in four variables:
\begin{equation}\label{eq4}
P_n(x,y,z,w)=\sum_{i,j}
P_n(i, j)x^i y^j z^{i+1} w^{n-2i-j},
\end{equation}
where $i$ and $j$ are in the same range as in \eqref{Pn}.

We also notice that this grammar can be used to investigate
the joint distribution of the number of peaks  and
the number of double descents of permutations on $[n]$.
To this end, we shall define polynomials $Q_n(x, y, z,w)$
and we shall show that the above grammar $G$ leads to a recurrence relation
on the polynomials $P_n(x,y)$, which also involves the polynomials
$Q_n(x,y,z,w)$.

Using the grammatical calculus,
 we deduce the following generating function of $P_n(x,y)$ without
 solving differential equations:
\begin{align}\label{gf1}
\sum_{n=0}^{\infty}P_n(x,y)\frac{t^n}{n!}=\frac{2\sqrt{(1 + y)^2 - 4x}
  e^{t/2\cdot(1-y+\sqrt{(1 + y)^2 - 4x})}}
  {1+y+\sqrt{(1 + y)^2 - 4x}-
  (1+y-\sqrt{(1 + y)^2 - 4x})e^{t \sqrt{(1 + y)^2 - 4x}}}.
\end{align}

We also establish a relationship between $\gen(z,t)$ and $\gen(y,t)$,
which leads to  the generating function of the polynomials $Q_n(x,y,z,w)$.  We show that
\begin{equation}\label{relation_1}
\gen(y,t)=y+xzF(x,y,z,w;t),
\end{equation}
where $F(x,y,z,w;t)$ is the generating function of the joint distribution of
the number of peaks, the number of valleys, the number of double descents and the number of double  rises, obtained by Carlitz and Scoville \cite{Carlitz-1974}. Thus we have given a grammatical
treatment of the formula of Carlitz and Scoville.

This paper is organized as follows.
In Section \ref{proof-of-Theorem-1.1}, we give an overview of
the formal derivative with respect to a context-free grammar and the
notion of a grammatical labeling.
Then a present a grammar $G$ and give a proof of Theorem \ref{thm1}
on the joint distribution of the number of exterior peaks and
the number of proper double descents.
 We also show that the same grammar $G$ also can be
 used to generate the polynomials $Q_n(x,y,z,w)$, which are the generating functions of the number of peaks and the numbers of double descents.
 Using the grammar, a recurrence relation of $P_n(x,y)$ is derived.
In Section 3, we give a proof of Theorem \ref{thm2} without solving differential
equations.


\section{A Grammatical Labeling of Permutations} \label{proof-of-Theorem-1.1}

In this section, we give an overview of the formal derivative with respect to a context-free grammar. Then we present a
context-free grammar and a grammatical labeling of permutations,
which can be employed to generate the number of permutations on $[n]$ with a given number of  exterior peaks and a given number of proper double descents.
Using a different grammatical labeling of permutations,
we show that the same grammar can be used to generate the number of permutations on $[n]$ with a given number of peaks and a given number of double descents. These two statistics have been studied by Carlitz and Scoville \cite{Carlitz-1974}.

A {\it context-free grammar} $G$ over a variable set $V$ is defined as a set of substitution rules replacing a variable in $V$ by a Laurent polynomial of variables in $V$. The polynomials
considered in this paper are assumed to be over real numbers.
The {\it formal derivative} $D$ with respect to $G$ is a linear operator
acting on Laurent polynomials in variables in $V$ such that
\[
D(uv)=uD(v)+vD(u).
\]
If $c$ is a constant,
we define $D(c) = 0$. Thus  we have $D(u^{-1}) = -u^{-2}D(u)$ since $D(uu^{-1}) = 0$. Clearly, we have the Leibniz formula:
\begin{equation}\label{leibniz}
  D^n(uv) =\sum_{k=0}^n {n\choose k}
D^k(u)D^{n-k}(v).
\end{equation}

For a Laurent polynomial $w$ of variables in $V$,
we define the generating function of  $w $ by
\[\gen(w,t)=\sum_{n=0}^\infty
D^n(w)\frac{t^n}{n!}.\]
Then the following relations hold:
\begin{align}
  \gen(u+v,t)&= \gen(u,t)+\gen(v,t),  \label{equ2}\\[6pt]
  \gen(uv,t) &= \gen(u, t)\gen(v, t),\label{equ1}\\[6pt]
  \gen'(u,t) &= \gen(D(u), t)\label{equ3},
\end{align}
where $u$, $v$ are Laurent polynomials of variables in $V$.

The idea of using the formal derivative with respect to
 a context-free grammar to study
combinatorial structures was initiated by Chen \cite{chen-1993}.
Dumont \cite{Dumont-1996} found grammars
for several classical combinatorial structures.
For example, the grammar
\begin{equation}\label{grammar1}
x \rightarrow xy,
\quad  y \rightarrow xy
\end{equation}
can be used to generate the Eulerian polynomials.
Moreover, Dumont \cite{Dumont-1996} gave the grammar
\[x\rightarrow xy, \quad  y \rightarrow x\]
to generate the Andr\'e polynomials $E_n(x,y)$.
Note that the generating function  for the Andr\'e polynomials $E_n(x,y)$ was first obtained by Foata and Sch\"utzenberger \cite{Foata-1973} by solving a differential equation.
However, Foata and Han \cite{Han-2001} later found a way to compute the generating function of $E_n(x, 1)$ without solving a differential equation.
Dumont \cite{Dumont-Ramamonjisop-1996} also discovered
the following grammar for the Ramanujan polynomials:
\[ x \rightarrow x^3y,
\quad  y \rightarrow xy^2.\]

Recently, the concept
of a {\it grammatical labeling}  was introduced in \cite{chen-2017}.
More precisely, a grammatical labeling is
 an assignment of the underlying elements of
a combinatorial structure with constants or variables,
 which is consistent with
 the substitution rules of a grammar $G$.
For example, given a context-free grammar
 \begin{equation}\label{grammar2}
 x \rightarrow xy,\quad y \rightarrow x^2,
 \end{equation}
we may label the elements of a permutation $\pi$ on $[n]$ by  $x$ and $y$ based on the exterior peaks.
Then it can be shown that $D^n(x)$ yields the generating function of the number of permutations on $[n]$ with a given number of exterior peaks.
Reminiscent of \eqref{grammar2},  Ma \cite{Ma-2012} found a connection between the number of peaks of permutations and the relations $D_z(x)=xy$ and $D_z(y)=x^2$, where $x=\sec(z)$, $y=\tan(z)$ and $D_z$ is the ordinary derivative with respect to $z$. Ma, Ma, Yeh and Zhu \cite{MaMaYehZhu-16} also found grammars to generate several polynomials associated with Eulerian polynomials, including $q$-Eulerian polynomials, alternating run polynomials and derangement polynomials.

Let $P_n(i, j)$ denote the number of  permutations on $[n]$ with $i$ exterior peaks and
$j$ proper double descents.
We now give a grammar to generate $P_n(i, j)$. Let
$P_n(x, y, z, w)$ be defined as in \eqref{eq4}, that is, \[P_n(x,y,z,w)=\sum_{i,j}
P_n(i, j)x^i y^jz^{i+1} w^{n-2i-j},\]
where $0\leq j \leq n-1$ and  $2i+j\leq n$.
The first few values of $P_n(x,y,z,w)$ are listed below:
 \begin{align*}
  P_1(x,y,z,w)&=zw,\\[6pt]
  P_2(x,y,z,w)&= z w^2 + x z^2,\\[6pt]
  P_3(x,y,z,w)&= z w^3+ 4 x z^2 w+ x y z^2,\\[6pt]
  P_4(x,y,z,w)&= z w^4 + 11 x z^2 w^2+ 6 x y z^2 w+ 5 x^2 z^3 + x y^2 z^2,\\[6pt]
  P_5(x,y,z,w)&= z w^5+ 26 x z^2 w^3+ 23 x y z^2 w^2+ 43 x^2 z^3 w+8 x y^2 z^2 w\\[6pt]
  &\,\,\qquad + 18 x^2 y z^3 + x y^3 z^2.
\end{align*}

Denote by $G$  the context-free grammar
\begin{align}
G\colon x \rightarrow xy,\quad y \rightarrow xz,\quad z \rightarrow zw,\quad w \rightarrow xz,
\end{align}
and let $D$ be the formal derivative with respect to  $G$.
This grammar can be viewed as a unification of grammars \eqref{grammar1} and \eqref{grammar2}.
Substituting $z, y$ by $x$ and substituting $x,w$  by $y$,
$G$ becomes the grammar \eqref{grammar1}.
Substituting $z, x$ by $x$, and substituting $w,y$  by $y$,
 $G$ becomes the grammar \eqref{grammar2}.


\begin{thm}\label{thm1}
For $n\geq 0$,
\begin{align}
  D^n(z)
=P_n(x, y,z, w).
\end{align}
\end{thm}

For example, for $n=4$, we have
\[
D^4(z)=z w^4 + 11 x z^2 w^2+ 6 x y z^2 w+ 5 x^2 z^3 + x y^2 z^2.
\]
The coefficient of $ x y z^2 w$
in $D^4(z)$ is 6,
corresponding to  the six permutations on $\{1, 2, 3, 4\}$
with one exterior peak and one proper double descent: 1432, 2431, 3214, 3421, 4213 and 4312.

To prove Theorem \ref{thm1}, we need a grammatical labeling  of a
permutation $\pi=\pi_1\pi_2\cdots \pi_n$ on $[n]$ by four variables $x, y,z, w$.
We first add an element $0$ at the end
of $\pi$ and label it by $z$.
If $i$ is an exterior peak, we label $\pi_i$ by $x$ and label $\pi_{i+1}$ by $z$,
if $i$ is a proper double descent,
we label $\pi_i$ by $y$. The rest of the elements
in $\pi$ are labeled by $w$.
Define the weight of $\pi$ by $$
w(\pi)=x^{\# {\rm exterior \, peaks }}y^{\# \, {\rm  proper\, double \,descents}}z^{\#  {\rm exterior \, peaks }+1}w^{\# \, {\rm otherwise}}.
$$
For example,  let $\pi=356412$. The labeling of $\pi$ is
as follows:
$$
\begin{array}{cllllll}
3 &5 &6 &4 &1 &2 &0\\
w &w &x &z &y & w & z
\end{array},
$$
and $w(\pi)=xyz^2w^3$.

\noindent{\it Proof of Theorem \ref{thm1}.} We proceed by induction on $n$. For $n=1$,
the grammatical labeling $\pi=1$ is given by
$$
\begin{array}{cl}
1 &0 \\
w &z
\end{array}.
$$
This yields $P_1(x,y,z,w)=wz$. On the other hand,
with respect to the grammar $G$, we have
$D(z)=zw$. Hence the theorem is valid for $n=1$.

 Assume that the theorem holds for $n$, that is, $D^n(z)=P_n(x,y,z,w)$. Let $\pi=\pi_1\pi_2\cdots \pi_n$ be a permutation on $[n]$ with $i$
 exterior peaks and $j$ proper double descents. Clearly, the weight $w(\pi)$
 is $x^iy^jz^{i+1}w^{n-2i-j}$.  We add a zero at the end of $\pi$. Then we insert $n+1$ into $\pi$ to generate a new permutation on $[n+1]$ with a zero at the end. According to where $n+1$ is inserted, there are four cases to label $n+1$ and relabel some elements in $\pi$.

\noindent Case 1:
 $n+1$ is inserted immediately before an exterior peak $k$ such that $\pi_k$ labeled by $x$ and $2\leq k\leq n-1$. Then we label $n+1$ by $x$ and relabel $\pi_k$ by $z$ and  $\pi_{k+1}$ by $y$. The change of labelings is illustrated as follows:
    \[\begin{array}{ccccc}
    \pi_{k-1}&<&\pi_k&>&\pi_{k+1}\\[-2pt]
    &&x&&z
    \end{array}
   \quad
   \Longrightarrow
   \quad
    \begin{array}{ccccccc}
    \pi_{k-1}&<&n+1&>&\pi_k&>&\pi_{k+1}\\[-2pt]
    &&x&&z&&y
    \end{array}.\]
    For example, if we insert $7$
before $6$ in the above example, we get
$$
\begin{array}{cccccccc}
3 &5 &7 &6 &4 &1 &2 &0\\
w &w &x &z &y &y & w & z
\end{array}.
$$
Similarly, for $i=1$ and $\pi_1>\pi_2$, the relabeling is
shown below:
\[\begin{array}{ccc}
    \pi_1&>&\pi_2\\[-2pt]
    x&&z
    \end{array}
   \quad
   \Longrightarrow
   \quad
    \begin{array}{ccccc}
    n+1&>&\pi_1&>&\pi_2\\[-2pt]
    x&&z&&y
    \end{array}.\]
Therefore, the insertion in this case always corresponds to the rule  $x \rightarrow x y$ and produces $i$ permutations on $[n+1]$ with $i$
exterior peaks and $j+1$ proper double descents. The sum of weight of these
$i$ permutations equals $ix^{i}y^{j+1}z^{i+1}w^{n-2i-j}$.

\noindent Case 2:  $n+1$ is inserted immediately before $\pi_k$ with labeling $y$, where $3\leq k\leq n$.  We label $n+1$ by $x$ and relabel $\pi_k$ by $z$, as illustrated below:
    \[\begin{array}{ccccc}
    \pi_{k-2}&>&\pi_{k-1}&>&\pi_k\\[-2pt]
    &&&&y
    \end{array}
   \quad
   \Longrightarrow
   \quad
    \begin{array}{ccccccc}
    \pi_{k-2}&>&\pi_{k-1}&<&n+1&>&\pi_k\\[-2pt]
    &&&&x&&z
    \end{array}.\]
For example, if we insert
$7$ before $2$ in the above example, we obtain
$$
\begin{array}{clllllll}
3 &5 &6  &4 &1 &7 &2  &0\\
w &w &x &z &y &x &z &z
\end{array}.
$$
Notice that this insertion corresponds to the rule $y \rightarrow xz$,
 and produces $j$ permutations on $[n+1]$ with $i+1$
exterior peaks and $j-1$ proper double descents. The sum of
weights of these permutations equals
 $jx^{i+1}y^{j-1}z^{i+2}w^{n-2i-j}$.

\noindent Case 3:  $n+1$ is inserted immediately before an element with  label $z$.  Here are two subcases. If $n+1$ is inserted before zero,
we just label $n+1$ by $w$:
    \[\begin{array}{cc}
    \pi_n&0\\[-2pt]
    &z
    \end{array}
   \quad
   \Longrightarrow
   \quad
    \begin{array}{cccc}
    \pi_n&<&n+1&0\\[-2pt]
    &&w&z
    \end{array}.\]
    For example, if we insert
$7$ before $0$, we get
$$
\begin{array}{clllllll}
3 &5 &6 &4 &1 &2 &7 &0\\
w &w &x &z &y & w & w  &z
\end{array}.
$$
If $n+1$ is inserted immediately before $\pi_k$, where $2\leq k\leq n$, we label $n+1$ by $x$ and relabel $\pi_{k-1}$ by $w$ as shown
below:
\[\begin{array}{ccc}
     \pi_{k-1}&>&\pi_k\\[-2pt]
     x&&z
    \end{array}
   \quad
   \Longrightarrow
   \quad
    \begin{array}{ccccc}
    \pi_{k-1}&<&n+1&>&\pi_k\\[-2pt]
    w&&x&&z
\end{array}.\]
For example, if we insert $7$ before $4$ in the above example, we get
$$
\begin{array}{clllllll}
3 &5 &6 &7 &4 &1  &2  &0\\
w &w &w &x &z &y &w  &z
\end{array}.
$$
In summary, the insertion in this case always corresponds to the rule $z \rightarrow zw$ and produces $i+1$ permutations on $[n+1]$ with $i$
exterior peaks and $j$ proper double descents. So we obtain a total weight   $(i+1)x^{i}y^{j}z^{i+1}w^{n-2i-j+1}$.

\noindent Case 4:  $n+1$ is inserted immediately before  $\pi_k$ with label $w$, where $1\leq k\leq n$. Then we label $n+1$ by $x$ and relabel $\pi_k$ by $z$:
    \[\begin{array}{c}
    \pi_{k}\\[-2pt]
    w
    \end{array}
   \quad
   \Longrightarrow
   \quad
    \begin{array}{ccc}
    n+1&>&\pi_k\\[-2pt]
    x&&z
    \end{array}.\]
For example, if we insert $7$
before $2$, we get
$$
\begin{array}{clllllll}
3 &7 &5 &6  &4 &1  &2  &0\\
w &x &z &x &z &y &w &z
\end{array}.
$$
This insertion corresponds to the rule $w \rightarrow x z$
and produces $n-2i-j$ permutations on $[n+1]$ with $i+1$ exterior peaks and $j$ proper double descents. So we get a total weight $(n-2i-j)x^{i+1}y^{j}z^{i+2}w^{n-2i-j-1}$.
Combining the above cases, we see that
\begin{align*}
 P_{n+1}(x,y,z,w)&=\sum_{i,j=0}^nP_n(i,j)\left(ix^{i}y^{j+1}z^{i+1}w^{n-2i-j}
+jx^{i+1}y^{j-1}z^{i+2}w^{n-2i-j}\right.\\[6pt]
&\left.\qquad\ +(i+1)x^{i}y^{j}z^{i}w^{n-2i-j+1}
+(n-2i-j)x^{i+1}y^{j}z^{i+2}w^{n-2i-j-1}\right).
\end{align*}
But, by the grammar $G$ we obtain that
\begin{align*}
D(x^iy^jz^{i+1}w^{n-2i-j})=&ix^{i}y^{j+1}z^{i+1}w^{n-2i-j}
+jx^{i+1}y^{j-1}z^{i+2}w^{n-2i-j}\\[6pt]
&\quad +(i+1)x^{i}y^{j}z^{i}w^{n-2i-j+1}
+(n-2i-j)x^{i+1}y^{j}z^{i+2}w^{n-2i-j-1}.
\end{align*}
Hence, by the induction hypothesis, we find that
\begin{align*}
D^{n+1}(z)&=D\left(P_n(x,y,z,w)\right)\\[6pt]
&=D\left(\sum_{i,j=0}^nP_n(i,j)x^iy^jz^{i+1}w^{n-2i-j}\right)\\[6pt]
&=\sum_{i,j=0}^nP_n(i,j)D(x^iy^jz^{i+1}w^{n-2i-j})\\[6pt]
&=\sum_{i,j=0}^nP_n(i,j)\left(ix^{i}y^{j+1}z^{i+1}w^{n-2i-j}
+jx^{i+1}y^{j-1}z^{i+2}w^{n-2i-j}\right.\\[6pt]
&\left.\qquad\ +(i+1)x^{i}y^{j}z^{i}w^{n-2i-j+1}
+(n-2i-j)x^{i+1}y^{j}z^{i+2}w^{n-2i-j-1}\right)\\[6pt]
&=P_{n+1}(x,y,z,w).
\end{align*}
Therefore, the theorem holds for $n+1$.
This completes the proof.
\qed

We note that
via a different grammatical labeling, the grammar
$G$ can also be used to deal with the
joint distribution of the number of peaks  and
the number of double descents of permutations on $[n]$.
 For a permutation
 $\pi=\pi_1\pi_2\cdots\pi_n$ on $[n]$, set $\pi_0=\pi_{n+1}=0$.
 For $1\leq i\leq n$, an index $i$ is called
 an  \emph{peak} if $\pi_{i-1}<\pi_i>\pi_{i+1}$, see \cite{Stanley-1997}.
 It is also called a \emph{maxima} in \cite{Carlitz-1974}, or a \emph{modified maximum}  in \cite{Goulden-2004}.
 An index $i$ is called a  \emph{double descent} or a \emph{double fall} if $\pi_{i-1}>\pi_i>\pi_{i+1}$, see \cite{Goulden-2004,Stanley-1997}.
For example, $\pi=4356721$ has two peaks: 1 and 5, and two double descents: 6 and 7.

 Let  $Q_n(i,j)$ denote the number of permutations on $[n]$ with $i$ peaks and $j$  double descents. Let
 \[Q_n(x, y,z, w)
=\sum_{i,j}
Q_n(i, j)x^iy^jz^{i} w^{n+1-2i-j},\]
where $0\leq j \leq n$ and  $2i+j\leq n+1$.
The first few values of $Q_n(x,y,z,w)$ are given below:
\begin{align*}
  Q_1(x,y,z,w)&=xz,\\[6pt]
  Q_2(x,y,z,w)&=xyz+xzw,\\[6pt]
  Q_3(x,y,z,w)&=xy^2z+2x^2z^2+2xyzw+xzw^2,\\[6pt]
  Q_4(x,y,z,w)&=xzw^3+3xyzw^2+8x^2z^2w+3xy^2zw+8 x^2yz^2+xy^3z.
\end{align*}

The following theorem shows that the
polynomials $Q_n(x,y,z,w)$ can also be generated by
the formal derivative $D$ with respect to the grammar $G$.

\begin{thm}\label{thm4}
 For $n\geq 1$,
\begin{align}
  D^n(y)
=Q_n(x, y,z, w).
\end{align}
\end{thm}

For example, it can be checked that
$D^4(y)=Q_4(x,y,z,w)$.
The coefficient  of $ x^2yz^2$ in $D^4(y)$ equals eight,
and there are eight  permutations on $\{1,2,3,4\}$ with
two peaks and one  double descent:
2143, 3142, 3214, 3241, 4132, 4213, 4231 and 4312.

The proof of the above theorem will be omitted.
Here we just provide a grammatical labeling.
Let $\pi=\pi_1\pi_2\cdots\pi_n$ be a permutation on $[n]$.
First, set $\pi_0=0$ and $\pi_{n+1}=0$.
For $1\leq i\leq n+1$,
if $i$ is a peak, then we label $\pi_i$ by $x$ and $\pi_{i+1}$ by $z$;
if $i$ is a  double descent,
then we label $\pi_{i+1}$ by $y$.
All other elements are labeled by $w$.
For example,  the permutation $\pi=4356721$ has the following labeling:
$$
\begin{array}{ccccccccc}
0 &4 &3 &5 &6 &7 &2  &1  &0\\
  &x &z &w &w &x &z  &y  &y
\end{array},
$$
and the weight of $\pi$ equals $x^2y^2z^2w^2$.

It should be noticed that Carliltz and Scoville defined the generating function on the joint distribution of the number of peaks, the number of valleys, the number of double descents and the number of double rises as:
 \[
 F_n(x,y,z,w)=\sum_{\pi}x^{\# {\rm peaks }-1}y^{\# \, {\rm double \,descents}}z^{\#  {\rm valleys}}w^{\# \, {\rm double \, rises}},
 \]
 where $\pi$ runs over the permutations on $[n]$.
 It turns out that the polynomials $Q_n(x,y,z,w)$ are essentially the
 polynomials $F_n(x,y,z,w)$ defined by Carlitz and Scoville. More precisely,
 it can be easily seen that for $n\geq 1$,
 \begin{equation}\label{relation_FQ}
 Q_n(x,y,z,w)=xzF_n(x,y,z,w).
 \end{equation}

The first few values of $F_n(x,y,z,w)$ are as follows:
\begin{align*}
  F_1(x,y,z,w)&=1,\\[6pt]
  F_2(x,y,z,w)&=y+w,\\[6pt]
  F_3(x,y,z,w)&=y^2+2xz+2yw+w^2,\\[6pt]
  F_4(x,y,z,w)&=w^3+3yw^2+8xzw+3y^2w+8 xyz+y^3.
\end{align*}

Let $F(x,y,z,w;t)$ denote the generating function of $F_n(x, y,z,w)$, that is,
\[F(x,y,z,w;t)
=\sum_{n=1}^\infty F_n(x,y,z,w)
\frac{t^n}{n!}.
\]
Carliltz and Scoville showed that
\begin{equation}\label{F_ex}
F(x,y,z,w;t)=
\frac{e^{vt}-e^{ut}}{ve^{ut}-ue^{vt}},
\end{equation}
where $uv=xz$ and $u+v=y+w$,
see  \cite[Exercise 3.3.46]{Goulden-2004} and \cite[Exercise 1.61]{Stanley-1997}.

 To conclude this section, we
 use the grammar $G$ to produce a recurrence relation of $P_n(x,y,z,w)$.

\begin{thm} For $n\geq 0$,
\begin{align}\label{P_{n+1}}
  P_{n+1}(x,y,z,w)
  =wP_{n}(x,y,z,w)+\sum_{k=0}^{n-1}{n\choose k}P_{k}(x,y,z,w)Q_{n-k}(x,y,z,w).
\end{align}
\end{thm}

\noindent {\it Proof.}  Since
\begin{equation}\label{D^n(zw)}
  D^{n+1}(z)=D^{n}(zw),
\end{equation}
by the Leibniz formula, we obtain that
\[D^{n+1}(z)=\sum_{k=0}^n{n\choose k}D^k(z)D^{n-k}(w).\]
But, for $k\geq 1$,
\begin{equation}\label{D^k(y)}
  D^k(w)=D^k(y),
\end{equation}
hence, by \eqref{D^n(zw)} and \eqref{D^k(y)}, we deduce that
\begin{align}\label{D^n+1(z)}
  D^{n+1}(z)=w D^{n}(z)+\sum_{k=0}^{n-1}{n\choose k}D^k(z)D^{n-k}(y).
\end{align}
Combining \eqref{D^n+1(z)} with Theorem \ref{thm1} and Theorem \ref{thm4},
we arrive at (\ref{P_{n+1}}). \qed

Setting $z=w=1$ in \eqref{P_{n+1}} and letting $Q_n(x,y)=Q_n(x,y,1,1)$,
we get
\begin{align}\label{eq6}
  P_{n+1}(x,y)
  =P_{n}(x,y)+\sum_{k=0}^{n-1}{n\choose k}P_{k}(x,y)Q_{n-k}(x,y).
\end{align}
Here we mention  some special cases of \eqref{eq6}.

Recall that $T_n(x)$ is defined by \eqref{thenumberofpeaks}.
Let $R(n,k)$ be the number of permutations on $[n]$
with $k$ peaks and let
\[R_{n}(x)=\sum_{k=0}^n R(n,k)x^k.\]
Taking $y=z=w=1$,  \eqref{P_{n+1}} becomes
\begin{equation}
T_{n+1}(x)=T_{n}(x)+\sum_{k=0}^{n-1}{n\choose k}T_{k}(x)R_{n-k}(x).
\end{equation}
Let $W(n,k)$ be the number of permutations on $[n]$
with $k$ double descents and let
\[W_{n}(y)=\sum_{k=0}^n W(n,k)y^k.\]
Taking $x=z=w=1$ in \eqref{P_{n+1}} yields
\begin{equation}
U_{n+1}(y)=U_{n}(y)+\sum_{k=0}^{n-1}{n\choose k}U_{k}(y)W_{n-k}(y),
\end{equation}
where $U_n(y)$ is defined by \eqref{yny}.

\section{The Generating Functions}\label{proof-of-Theorem-1.2}

In this section, we use the grammar
\[
G\colon \,\,  x \rightarrow xy,\quad y \rightarrow x z, \quad z \rightarrow z w, \quad  w \rightarrow  xz
\]
to derive the generating function of $P_n(x,y)$ without solving differential equations.
In fact, we shall consider the generating function $\gen(z,t)$ of
the polynomials $D^n(z)$ in four variables $x,y,z,w$, that is,
\[\gen(z,t)=\sum_{n=0}^{\infty}D^n(z)\frac{t^n}{n!}.\]
Furthermore, we show that the generating function
for the joint distribution of the number of peaks and the number of double descents can also be determined by $\gen(y,t)$. This leads to a grammatical approach
to the generating function of Carliltz and Scoville.

\begin{thm}\label{thm2}
  We have
  \begin{equation*}
  \gen(z,t)=\frac{2z\sqrt{(w + y)^2 - 4xz}
  e^{t/2\cdot(w-y+\sqrt{(w + y)^2 - 4xz})}}
  {w+y+\sqrt{(w + y)^2 - 4xz}-
  (w+y-\sqrt{(w + y)^2 - 4xz})e^{t\sqrt{(w + y)^2 - 4xz}}}.
  \end{equation*}
\end{thm}

Combining Theorem \ref{thm1} and Theorem \ref{thm2}, we readily deduce
 \eqref{gf1} by setting $w=z=1$.

To prove Theorem \ref{thm2}, we need the following relation
 between $\gen(z,t)$ and $\gen(x,t)$.

\begin{lem}\label{thm3} We have
\begin{equation}\label{eqnew1}
{\rm Gen}(z,t)=zx^{-1} {\rm Gen}(x,t)e^{(w-y)t}.
\end{equation}
\end{lem}

\noindent{\it Proof.} For the grammar $G$ and its
formal derivative $D$, it is evident that
\begin{align}\label{D(w-y)=0}
  D(w-y)=0.
\end{align}
By the Leibniz rule \eqref{leibniz}, we get
\begin{align*}
  D(x^{-1})=-x^{-2}D(x)=-x^{-2}xy=-x^{-1}y.
\end{align*}
It follows that
\begin{align}\label{w-y}
 D(zx^{-1})=zx^{-1}(w-y).
\end{align}
Combining \eqref{D(w-y)=0} and \eqref{w-y}, we deduce that
$$D^n(zx^{-1})=zx^{-1}(w-y)^n.$$
Thus,
\begin{equation}\label{gen(zx^-1)}
\gen(zx^{-1},t)=\sum_{n=0}^{\infty}D^n(zx^{-1})\frac{t^n}{n!}=zx^{-1}\sum_{n=0}^{\infty}(w-y)^n\frac{t^n}{n!}=zx^{-1}e^{(w-y)t}.
\end{equation}
By \eqref{equ1}, we see that
\begin{equation}\label{gen(zx-1x)}
\gen(z,t)=\gen(z x^{-1},t) \gen(x,t).
\end{equation}
Substituting \eqref{gen(zx^-1)} into \eqref{gen(zx-1x)}, we find that
\[
\gen(z,t)=zx^{-1}\gen(x,t)e^{(w-y)t}.
\]
This completes the proof.
\qed

\noindent{\it Proof of Theorem \ref{thm2}.}
By \eqref{equ1}, we see that
\begin{equation}\label{Inequation}
\gen(z x,t)={1\over \gen(z^{-1}x^{-1},t)}.
\end{equation}
Applying Lemma \ref{thm3}, we find that
\begin{equation}\label{gen(zx)}
\gen(zx,t)=\gen(z,t)\gen(x,t)=xz^{-1} e^{(y-w)t}\gen^2(z,t).
\end{equation}
Combining \eqref{Inequation} and \eqref{gen(zx)}, we deduce that
\begin{equation}\label{eq10}
\gen(z,t)=\sqrt{\frac{zx^{-1} e^{(w-y)t}}{{\rm Gen}(z^{-1}x^{-1},t)}}.
\end{equation}
We continue to compute $\gen(x^{-1}z^{-1},t)$.
It is easily checked that
\begin{align}
  D(x^{-1}z^{-1})&=-x^{-1}z^{-1}(w+y),\nonumber \\[6pt]\nonumber
  D^2(x^{-1}z^{-1})&=x^{-1}z^{-1}((w+y)^2-2xz),\\[6pt] \label{d3}
  D^3(x^{-1}z^{-1})&=-x^{-1}z^{-1}(w+y)((w+y)^2-4xz),\\[6pt] \label{d4}
  D^4(x^{-1}z^{-1})&=-x^{-1}z^{-1}(w+y)((w+y)^2-2xz)((w+y)^2-4xz).
\end{align}
Observe that
\begin{align}
  D((w+y)^2-4xz)=0.\label{fact1}
\end{align}
In light of \eqref{fact1}, it follows from
 (\ref{d3}) and (\ref{d4})  that for $k \geq 0$,
\begin{equation*}
 \begin{split}
 D^{2k+1}(x^{-1}z^{-1})&=-x^{-1}z^{-1}(w+y)((w+y)^2-4xz)^k, \\[6pt]
 D^{2k+2}(x^{-1}z^{-1})&=x^{-1}z^{-1}((w+y)^2-2x z)((w+y)^2-4xz)^k.
 \end{split}
\end{equation*}
Thus,
\begin{align*}
{\rm Gen}(x^{-1}z^{-1},t)&= \sum_{k=0}^\infty \frac{D^k(z^{-1}x^{-1})t^k}{k!}=\frac{1}{xz}\left(-\frac{2xz}{(w+y)^2-4xz}\right.\\[6pt]
&\qquad +\left(\frac{w+y}{2\sqrt{(w+y)^2-4xz}}+\frac{(w+y)^2-2xz}{2((w+y)^2-4xz)}\right)e^{-t\sqrt{(w+y)^2-4xz}}\\[6pt]
&\qquad  \left.-\left(\frac{w+y}{2\sqrt{(w+y)^2-4xz}}-\frac{(w+y)^2-2xz}{2((w+y)^2-4xz)}\right)e^{t\sqrt{(w+y)^2-4xz}}\right).
\end{align*}
Putting ${\rm Gen}(x^{-1}z^{-1},t)$ into \eqref{eq10} gives the required  formula for  $\gen(z,t)$. This completes the proof.
\qed

It can be seen that Theorem \ref{thm2} serves as a unification of
  Theorems \ref{lem1}, \ref{lem2} and \ref{lem3}. Taking $y=z=w=1$,
   Theorem \ref{thm2} reduces to Theorem \ref{lem1}. Setting  $x=z=w=1$
    in Theorem \ref{thm2} yields Theorem \ref{lem2}. Theorem
     \ref{thm2} simplifies to Theorem \ref{lem3} when $x=z=w=1$ and $y=0$.

We note that using the grammar $G$, it is easy to
establish a connection between the
generating functions of $P_n(x,y,z,w)$ and $Q_n(x,y,z,w)$.

\begin{thm}\label{relation_thm}
  We have
  \begin{equation}\label{relation_eqt}
    \gen(y,t)=\ln'(\gen(z,t))-w+y.
  \end{equation}
\end{thm}

\pf From the substitution rule $z\rightarrow zw$, we get
\[
\gen'(z,t)=\gen(D(z),t)=\gen(zw,t)=\gen(z,t)\gen(w,t),
\]
and hence
\[
  \gen(w,t)=\frac{\gen'(z,t)}{\gen(z,t)}=\ln'(\gen(z,t)).
\]
Since $D(w)=D(y)$, we see that
\[
\gen(y,t)=\gen(w,t)-w+y=\ln'(\gen(z,t))-w+y.
\]
This completes the proof. \qed

In view of Theorem \ref{relation_thm}, from the
generating function of $P_n(x,y,z,w)$ we can deduce the generating function of $Q_n(x,y,z,w)$.

\begin{thm}\label{coro}
  We have
  \begin{equation}\label{geny_eqt}
  \gen(y,t)=y+\frac{-2xz+2xze^{t\sqrt{(y+w)^2-4xz}}}
  {w+y+\sqrt{(w + y)^2 - 4xz}-
  (w+y-\sqrt{(w + y)^2 - 4xz})e^{t\sqrt{(w + y)^2 - 4xz}}}.
  \end{equation}
\end{thm}
\pf
By Theorem \ref{thm2}, we find that
\begin{align*}
  \ln(\gen(z,t))&=\ln(2z)+\frac{1}{2}\ln((w+y)^2-4xz)+\frac{t}{2}\left(w-y+\sqrt{(w+y)^2-4xz}\right)\\[4pt]
                &\qquad -\ln\left(w+y+\sqrt{(w + y)^2 - 4xz}\right.\\[4pt]
                &\qquad \qquad \qquad -
  \left.(w+y-\sqrt{(w + y)^2 - 4xz})e^{t\sqrt{(w + y)^2 - 4xz}}\right).
\end{align*}
Hence
\begin{equation}\label{lngenz_eqt}
\begin{split}
  \ln'(\gen(z,t))&=\frac{1}{2}\left(w-y+\sqrt{(w+y)^2-4xz}\right)\\
                &\qquad +\frac{(w+y-\sqrt{(w + y)^2 - 4xz})\sqrt{(w + y)^2 - 4xz}e^{t\sqrt{(w + y)^2 - 4xz}}}{w+y+\sqrt{(w + y)^2 - 4xz}-
  (w+y-\sqrt{(w + y)^2 - 4xz})e^{t\sqrt{(w + y)^2 - 4xz}}}.
\end{split}
\end{equation}
Plugging \eqref{lngenz_eqt} into \eqref{relation_eqt}, we obtain
\begin{align*}
  \gen(y,t)&=y-\frac{1}{2}\left(y+w-\sqrt{(w+y)^2-4xz}\right)\\
                &\qquad +\frac{(w+y-\sqrt{(w + y)^2 - 4xz})\sqrt{(w + y)^2 - 4xz}e^{t\sqrt{(w + y)^2 - 4xz}}}{w+y+\sqrt{(w + y)^2 - 4xz}-
  (w+y-\sqrt{(w + y)^2 - 4xz})e^{t\sqrt{(w + y)^2 - 4xz}}},
\end{align*}
which yields \eqref{geny}. This completes the proof. \qed

To conclude this paper, we remark that the above formula for
the generating function of $Q_n(x, y,z,w)$ can be
recast as the formula \eqref{F_ex}
of Carlitz and Scoville
for the generating function of $F_n(x,y,z,w)$.

Using \eqref{relation_FQ}, that is
\[  Q_n(x,y,z,w)=xzF_n(x,y,z,w),\]
we get
\begin{equation}\label{relation}
\gen(y,t)=y+xzF(x,y,z,w;t).
\end{equation}
Assume that $u+v=y+w$ and $uv=xz$, as in \eqref{F_ex}.
Then \eqref{geny_eqt} can be rewritten as
\begin{equation}\label{geny}
\begin{split}
\gen(y,t)&=y+\frac{uve^{t(u-v)}-uv}{u-ve^{t(u-v)}}\\[6pt]
            &=y+uv\frac{e^{tu}-e^{tv}}{ue^{tu}-ve^{tu}}.
\end{split}
\end{equation}
Thus \eqref{F_ex} follows from \eqref{geny} and \eqref{relation}.


\begin{thebibliography}{999}
  \bibitem{Aguiar-2004}
   M. Aguiar, N. Bergeron and K. Nyman, The peak algebra and the descent algebra of type $B$ and $D$, Trans. Amer. Math. Soc. 356 (2004), 2781--2824.

  \bibitem{Barry-2014}
  P. Barry, Constructing exponential Riordan arrays from their $A$ and $Z$ sequences, J. Integer Seq. 17 (2014), 14.2.6.

  \bibitem{Basset-2014}
  N. Basset, Counting and generating permutations using timed languages, Latin American Symposium on Theoretical Informatics, Springer, Berlin, Heidelberg, 2014, 502--513.


  \bibitem{Carlitz-1974}
  L. Carlitz and R. Scoville, Generalized Eulerian numbers: combinatorial applications, J. Reine Angew. Math. 265 (1974), 110--137.

  \bibitem{chen-1993}
  W.Y.C. Chen, Context-free grammars, differential operators and formal power series, Theoret. Comput. Sci. 117 (1993), 113--129.

  \bibitem{chen-2017}
  W.Y.C. Chen and A.M. Fu, Eulerian grammars, Permutations and increasing trees,
  Adv. in Appl. Math. 82 (2017), 58--82.

  \bibitem{David-1962}
  F.N. David and D.E. Barton, Combinatorial Chance, Hafner, New York, 1962.

  \bibitem{Dumont-1996}
   D. Dumont, Grammaires de William Chen et d\'erivations dans les arbres et arborescences, S\'em. Lothar. Combin. 37 (1996) B37a, 21 pp.
  \bibitem{Dumont-Ramamonjisop-1996}
   D. Dumont and A. Ramamonjisoa, Grammaire de Ramanujan et Arbres de Cayley, Electron. J. Combin. 3 (2) (1996), \#R17.

  \bibitem{Elizalde-2003}
  S. Elizalde and M. Noy, Consecutive patterns in permutations, Adv. in Appl. Math. 30 (2003), 110--125.

  \bibitem{Han-2001}
  D. Foata and G.-N. Han, Arbres minmax et polyn\^omes d'Andr\'e, Adv. in Appl. Math. 27 (2001) 367--389.

  \bibitem{Foata-1973}
    D. Foata and M.-P. Sch\"utzenberger, Nombres d'Euler et permutations alternantes, in: J.N. Srivastava, et al. (Eds.), A Survey of Combinatorial Theory, North-Holland, Amsterdam, 1973, pp. 173--187.

  \bibitem{Gessel-1977}
  I.M. Gessel, Generating functions and enumeration of sequences, Ph.D. Dissertation, Massachusetts Institute of Technology, 1977.

  \bibitem{Goulden-2004}
  I.P. Goulden and D.M. Jackson, Combinatorial Enumeration, Courier Corporation, 2004.

  \bibitem{Kuznetsov-1994}
  A.G. Kuznetsov, I.M. Pak and A.E. Postnikov, Increasing trees and alternating permutations, Russian Math. Surveys 49 (1994), 79--114.

  \bibitem{Ma-2012}
  S.-M. Ma, Derivative polynomials and enumeration of permutations by number of interior and left peaks, Discrete Math. 312 (2012), 405--412.

  \bibitem{MaMaYehZhu-16}
  S.-M. Ma, J. Ma, Y.-N. Yeh and B.-X Zhu, Context-free grammars for several polynomials associated with Eulerian polynomials, arXiv:1609.05829.


  \bibitem{Stanley-1997}
  R.P. Stanley, Enumerative Combinatorics, Vol. 1, Cambridge Stuies in Advanced Mathematics, 2012.

  \bibitem{Gessel-OEIS}
  The On-Line Encyclopedia of Integer Sequences, \url{https://oeis.org/A008971}.
\end{thebibliography}
\end{document}